\documentclass{article}
\usepackage{amsmath,array}
\usepackage{amsfonts}
\usepackage{hyperref}
\usepackage{epsf}
\usepackage[dvips]{graphicx}
\usepackage{wrapfig}
\usepackage{caption}
\makeatletter

\newtheorem{thm}{Theorem}
\newtheorem{prp}{Proposition}[section]
\newtheorem{rmk}{Remark}
\@addtoreset{equation}{section}
\makeatother
\begin{document}

\begin{center}
{\Large \bf { Invariant Subspace Method and Fractional\\ Modified Kuramoto-Sivashinsky Equation}}\\
\vskip1.4 truecm {\bf  A. Ouhadan$^{1}$ and  E. H. El Kinani$^{2,\dagger}$}

\vspace{0.5cm}

$^{1}$ Centre R\'egional des M\'etiers de l'Education et de la Formation,  Mekn\`es, \\BP 255, Morocco.\\

$^{2}$A.A Group, Mathematical Department Moulay Isma\"{\i}l University, Faculty of Sciences and Technics Errachidia, BP 509, Morocco.\\
$^{\dagger}$ Universit\'e Moulay Isma\"{\i}l Ecole Nationale Sup\'erieure des Arts et M\'etiers (ENSAM), Marjane 2, B.P. 15290, Mekn\`es, Maroc.\\

\vspace{0.5cm}

\end{center}

\begin{abstract}
In this paper, the  invariant subspace method is applied to the time fractional modified Kuramoto-Sivashinsky partial differential equation. The obtained reduced system of nonlinear ordinary fractional equations is solved by the Laplace transform method and with using of some useful properties of Mittag-Leffler function. Then, some exact solutions of the time fractional nonlinear  studied equation are found.\end{abstract}

\vspace{1cm}

{\bf{Keywords}}\\

Invariance subspace method, Caputo fractional derivative,  fractional  modified Kuramoto-Sivashinsky equation,  Mittag-Lefller function.
\begin{section}{Introduction}
In the last decade, fractional calculus  attracted a great interest of many researchers. The idea of fractional order derivative was started with half-order derivative as discussed  in the literature by Leibniz and L'H$\hat{o}$pital. Next, it was extended to an arbitrary order derivative by  Liouville, Riemann, Gr$\ddot{u}$nwald, Letnikov, Caputo etc. In addition,  different approaches to define fractional derivatives are  known \cite{1,2,3,4}. The study of fractional differential equations becomes of great interest, since for their widely applications including fluid flow, dynamical
processes in self-similar and porous structures,  electromagnetic waves, probability and statistics,  viscoelasticity,  signal processing, and so on \cite{1,4,5}.\\

The construction of particular exact solutions of  fractional differential equations is not an easy task and it remains a relevant problem. This is the reason why a  powerful methods for solving those fractional  equations were recently developed in the literature, including
Adomian decomposition method \cite{6}, first integral method \cite{7}, homotopy perturbation method \cite{8}, Lie group theory method \cite{9,10,11,12,13} and so on. Most recently, according to invariance principles, the invariant subspace method developed by V.A. Galaktionov and S.R. Svirshchevski \cite{14} to study partial differential equations was extended by R.K. Gazizov and A.A. Kasatkin \cite{15} to construct some particular exact solutions for time fractional differential equations.\\

The invariant subspace method used in the present paper, yields us with  an exact  solutions of the time fractional modified Kuramoto-Sivashinsky equation in terms of the well known Mittag-Leffler function. In the paper \cite{15}, the invariant subspace method and Lie group analysis are joined to solve the reduced fractional ordinary differential system and the original studied equation. In our case, resolution of the reduced system is done by the Laplace transform method and  by using  of some remarkable properties of the well known Mittag-Leffler function \cite{16,17,18}.\\

This paper is organized as follows: In section 2, we recall some main results of fractional derivatives and integrals. Section 3, is devoted to describe the invariance subspace method. While in  section 4, we use the described method to construct exact solution admitted by the  time fractional modified Kuramoto-Sivashinsky equation. Finally, a conclusion is given.
\end{section}
\begin{section}{SOME BASIC RESULTS ON FRACTIONAL CALCULUS}
This section is devoted to recall briefly some definitions and basic results on fractional calculus. For more details and proofs of the results, we refer to \cite{1,2,3,4}.\\

The Riemann-Liouville fractional integral is defined by:
\begin{equation}J^{\gamma}_{t}f(t)=\frac{1}{\Gamma(\gamma)}\int_{0}^{t}(t-\tau)^{\gamma-1}f(\tau)\,\,d\tau,\label{integrl}\end{equation}
where  $\gamma\in\mathbb{R^{+}}$, and \begin{equation}\Gamma(\gamma)=\int_{0}^{+\infty}x^{\gamma-1}e^{-x}\,dx,\end{equation}
is the Euler Gamma function.\\

By definition $J^{0}_{t}f(t)=f(t)$ and it satisfies the property $J^{\alpha}_{t}J^{\gamma}_{t}f(t)=J^{\alpha+\gamma}_{t}f(t)$.\\

 First recall that there are various contributions  \cite{1,2,3,4} to define  fractional derivatives. In this paper, we adopt the fractional  derivative in the sense of Caputo \cite{1,2,3,4}. The Caputo definition is used not only because it makes easy the consideration of initial conditions but also because the derivative of a constant is equal to zero. In what follows, we recall some   important results and properties  of fractional derivatives and integrals. For more details see for example \cite{4}. Before going on, let us  denote by $AC^{n}([0,t]), n\in\mathbb{N}$ the class  of functions $f(x)$ which are continuously differentiable in $[0,t]$ up to order $(n-1)$ and with $f^{(n-1)}\in AC([0,t])$.\\
\begin{thm}
 Let $n-1<\alpha<n$, with $n\in\mathbb{N}$. If $f(x)\in AC^{n}([0,t])$, then the Caputo fractional derivative exists almost everywhere on $[0,t]$ and it is represented in the form:
 \begin{equation}D^{\alpha}_{t}f(t)=J^{n-\alpha}_{t}D^{n}_{t}f(t)=\frac{1}{\Gamma(n-\alpha)}\int_{0}^{t}(t-\tau)^{n-\alpha-1}f^{(n)}(\tau)\,d\tau,\quad \alpha\neq n. \label{caputodef}\end{equation}
\end{thm}
The Caputo derivative (\ref{caputodef}) and the Riemann-Liouville integral (\ref{integrl})satisfy the following properties \cite{4}:
 \begin{eqnarray}
 D^{\alpha}_{t}J^{\alpha}_{t}f(t)&=&f(t),\qquad \alpha>0,\\
 J^{\alpha}_{t}D^{\alpha}_{t}f(t)&=&f(t)-\sum_{k=0}^{n-1}f^{(k)}(0)\frac{t^{k}}{k!},\qquad \alpha>0,\quad t>0,\label{integderiv}\\
 J^{\alpha}_{t}t^{\gamma}&=&\frac{\Gamma(\gamma+1)}{\Gamma(\gamma+\alpha+1)}t^{\gamma+\alpha},\qquad \alpha>0,\quad\gamma>-1,\quad t>0,\\
 D^{\alpha}_{t}t^{\gamma}&=&\frac{\Gamma(\gamma+1)}{\Gamma(\gamma-\alpha+1)}t^{\gamma-\alpha}, \alpha>0,\gamma\in]-1,0[\cup]0,+\infty[,\, t>0.
 \end{eqnarray}

Here, it is important to mention that the studied equation is a time fractional partial differential equation of order $0<\alpha<1$, so the integer $n$ appearing in the relation (\ref{caputodef}) is equal to one. Consequently, the formulae (\ref{integderiv}) becomes:
\begin{equation}J^{\alpha}_{t}D^{\alpha}_{t}f(t)=f(t)-f(0),\qquad  t>0.\end{equation}
\end{section}
\begin{section}{DESCRIPTION OF THE INVARIANT SUBSPACE  METHOD}

The aim of this section is to collect and to present some necessary and essentials results from invariant subspace theory.
The invariant subspace method \cite{14} was firstly used to construct particular exact solutions of evolutionary partial differential equations of the form:
\begin{equation}\frac{\partial u}{\partial t}=F(u,u_{1x}, u_{2x},\ldots, u_{kx}),\quad k\in\mathbb{N},\end{equation}where $u=u(t,x), u_{ix}=\frac{\partial^{i} u}{\partial x^{i}}$  is the i-th order  derivative of $u$ with respect to the space variable $x$ and $F$ is a nonlinear differential operator.\\

Recently, Gazizov and Kasatkin \cite{15} showed that the invariant subspace method can be applied also  to equations with time fractional derivative:\\

In fact, consider the time fractional partial differential equation of the form:
\begin{equation}D^{\alpha}_{t}u(t,x)=F[u],\quad \label{Eqg1}\end{equation}
where $F[u]=F(u,u_{1x}, u_{2x},\ldots, u_{kx})$ and $D^{\alpha}_{t}$ is the time  fractional derivative in the sense of Caputo. The invariant subspace method is based on the following basic definitions and results \cite{14, 15}.\\

Let $f_{1}(x),\ldots,f_{n}(x)$ be an  $n$ linearly independent functions and $W_{n}$ is the $n$-dimensional linear space namely $W_{n}=\langle f_{1}(x),\ldots,f_{n}(x)\rangle$. $W_{n}$ is said to be invariant under the given operator $F[u]$ if $F[u]\in W_{n}$ whenever $u\in W_{n}$.

\begin{prp} Let $W_{n}$ be  an invariant subspace of $F[u]$. A function $u(t,x)=\sum_{i=1}^{n}f_{i}(x)u_{i}(t)$ is a solution of equation (\ref{Eqg1}) if and only if the expansion coefficients $u_{i}(t)$ satisfy the following system of fractional ordinary differential equations:
\begin{equation}\left\{\begin{array}{ccc} D^{\alpha}_{t}u_{1}&=&F_{1}(u_{1},\ldots,u_{n}),\nonumber\\
D^{\alpha}_{t}u_{2}&=&F_{2}(u_{1},\ldots,u_{n}),\nonumber\\
\vdots&&\vdots\nonumber\\
D^{\alpha}_{t}u_{n})&=&F_{n}(u_{1},\ldots,u_{n}),\end{array}\right.\end{equation}
where $F_{1},\ldots,F_{n}$ are given by: 
\begin{equation}
F(c_{1}f_{1}(x)+\cdots+c_{n}f_{n}(x))=F_{1}(c_{1},\ldots,c_{n})f_{1}(x)+\cdots+F_{n}(c_{1},\ldots,c_{n})f_{n}(x).
\end{equation}
\end{prp}
\begin{rmk}A crucial question in the theory of invariant subspace method was how to get the corresponding invariant subspace of a given differential operator. This question is solved by the following proposition and for more details see \cite{14}.\end{rmk}

\begin{prp} Let $f_{1}(x),\ldots,f_{n}(x)$ form the fundamental set of solutions of a linear n-th order ordinary differential equation
\begin{equation}L[u]=y^{(n)}+a_{1}(x)y^{(n-1)}+\cdots+a_{n-1}(x)y'+a_{n}(x)y=0,\quad \label{Eqdl1}\end{equation}
and $F[y]=F(x,y,y',\ldots,y^{(k)})$ a given differential operator of order $k\leq n-1$, then the subspace $W_{n}=\langle f_{1}(x),\ldots,f_{n}(x)\rangle$ is invariant with respect to $F$ if and only if: \begin{equation}L[F[y]]=0,\end{equation}
whenever $y$ satisfies equation (\ref{Eqdl1}).\end{prp}

\begin{rmk}Condition of invariance appearing in the above proposition is  the invariance criterion for equation (\ref{Eqdl1})
with respect to the Lie-B$\ddot{a}$cklund generator $V=F[y]\frac{\partial}{\partial y}$. This criterion shows us how the invariant subspace method is related  to the techniques used in Lie symmetry analysis, see for more details \cite{19,20,21,22}.\end{rmk}
\end{section}
\begin{section}{EXACT SOLUTION OF THE FRACTIONAL MKS EQUAION}
In this section, we use the invariant subspace method to construct some exact solutions of the  time fractional modified Kuramoto-Sivashinsky equation (mKS) which is given by:
\begin{equation}D^{\alpha}_{t}(u)=-u_{4x}-u_{2x}+(1-\lambda)(u_{x})^{2}+\lambda(u_{xx})^{2},\quad0<\alpha<1,\, t>0, \label{mkseq}\end{equation}
where $u=u(t,x)$ and $\lambda\in]0;1[$. In the case $\alpha=1$, the (mKS) equation (\ref{mkseq}) is a model for the dynamics of a
hyper-cooled melt \cite{23,24}. A more general class of such models was introduced and discussed
in \cite{25}.

\begin{prp}For any $\lambda\in]0;1[$ the  nonlinear operator $F[u]$  given  by:
\begin{equation}F[u]=-u_{4x}-u_{2x}+(1-\lambda)(u_{x})^{2}+\lambda(u_{xx})^{2},\end{equation}admits $W_{3}=\langle 1,\cos\gamma x,\sin\gamma  x\rangle$ with $\gamma=\sqrt{\frac{1-\lambda}{\lambda}}$ as an invariant subspace.\end{prp}

\noindent{}\textbf{Proof.} For any function  \begin{equation}h(t,x)=C_{1}+C_{2}\cos\gamma x+C_{3}\sin\gamma  x,\end{equation}
with $C_{i}=C_{i}(t)$ arbitrary functions, we get:
\begin{eqnarray}F[h]&=&-\gamma^{4}C_{2}\cos\gamma x-\gamma^{4}C_{3}\sin\gamma x+\gamma^{2}C_{2}\cos\gamma x+\gamma^{2}C_{3}\sin\gamma x\nonumber\\
&&+(1-\lambda)(\gamma^{2}C_{3}^{2}\cos^{2}\gamma x+\gamma^{2}C_{2}^{2}\sin^{2}\gamma x-2\gamma^{2}C_{2}C_{3}\cos\gamma x\sin\gamma x)\nonumber\\
&&+\lambda(\gamma^{4}C_{2}^{2}\cos^{2}\gamma x+\gamma^{4}C_{3}^{2}\sin^{2}\gamma x+2\gamma^{4}C_{2}C_{3}\cos\gamma x\sin\gamma x)\nonumber\\
&=&(1-\lambda)\gamma^{2}(C_{2}^{2}+C_{3}^{2})+(\gamma^{2}-\gamma^{4})C_{2}\cos\gamma x+(\gamma^{2}-\gamma^{4})C_{3}\sin\gamma x\in W_{3}.\nonumber
\end{eqnarray}
Now, we search  an exact solution   admitted by the time fractional (mKS) equation (\ref{mkseq}) of the form:
\begin{equation}u(t,x)=C_{1}+C_{2}\cos\gamma x+C_{3}\sin\gamma  x.\quad \label{formsolut}\end{equation}
Consequently, a function $u(t,x)$ of the form (\ref{formsolut}) is a solution of the time fractional (mKS) equation if the expansion coefficients $C_{i}(t)$ satisfy the following system of ordinary fractional differential equations:
\begin{equation}\left\{\begin{array}{ccl} D^{\alpha}_{t}C_{1}&=&(1-\lambda)\gamma^{2}C_{3}^{2}+(1-\lambda)\gamma^{2}C_{2}^{2},\\
D^{\alpha}_{t}C_{2}&=&\gamma^{2}(1-\gamma^{2})C_{2},\\
D^{\alpha}_{t}C_{3}&=&\gamma^{2}(1-\gamma^{2})C_{3}.\end{array}\right.\label{system}\end{equation}
To get a non trivial solution needs to assume the condition $C_{2}(0)C_{3}(0) \neq 0$ and for convenience we suppose $C_{2}(0)=C_{3}(0)=1$. This last condition will be clear when the Laplace transform will be used. We start to  construct  solution of the third equation in the  above reduced system of ordinary fractional differential equations. We mention that, with the Laplace transform it is frequently possible to avoid
working with equations of different differential orders  by translating the problem into an easy one.\\

Recalling  some useful properties of the Laplace transform \cite{1}:
\begin{equation}\mathfrak{L}\left\{D^{\alpha}_{t}f(t)\right\}=s^{\alpha}\widetilde{f}(s)-s^{\alpha-1}f(0),\quad 0<\alpha<1, \label{laplaceprop1}\end{equation}
where \begin{equation}\mathfrak{L}\left\{f(t)\right\}=\widetilde{f}(s)=\int_{0}^{\infty}e^{-st}f(t)dt.\end{equation}
By putting $\theta=\gamma^{2}(1-\gamma^{2})$ and applying the Laplace transform on both sides of the third equation appearing in the fractional ordinary differential system, we obtain:
\begin{equation}s^{\alpha}\mathfrak{L}\left\{C_{3}(t)\right\}-s^{\alpha-1}C_{3}(0)=\theta \mathfrak{L}\left\{C_{3}(t)\right\},\end{equation}it yields:
\begin{equation}\mathfrak{L}\left\{C_{3}(t)\right\}(s)=\frac{s^{\alpha-1}}{s^{\alpha}-\theta},\end{equation}
then, with the inverse Laplace transform, it gives:
\begin{equation}C_{3}(t)=E_{\alpha,1}(\theta t^{\alpha}),\quad \label{expressionC3}\end{equation}
where $E_{\alpha,\beta}(\theta t^{\alpha})$ is the Mittag-Leffler function given by:
\begin{equation}E_{\alpha,\beta}(\theta t^{\alpha})=\Sigma_{i=0}^{\infty}\frac{\theta^{i}t^{\alpha i}}{\Gamma(\alpha i+\beta)}.\quad\end{equation}
Not that when $\beta=1$, $E_{\alpha,1}\equiv E_{\alpha}$.\\

Two last equations in the fractional ordinary differential system (\ref{system}) are the same, hence,
 \begin{equation}C_{2}(t)=C_{3}(t)=E_{\alpha,1}(\theta t^{\alpha}).\end{equation}
Substituting the obtained expressions of  $C_{2}$ and $C_{3}$ in the first equation of the system (\ref{system}), it leads to:
\begin{equation}D^{\alpha}_{t}C_{1}=2(1-\lambda)\gamma^{2}\left( E_{\alpha,1}(\theta t^{\alpha})\right)^{2}.\label{equadeC1}\end{equation}
The Mittag-Leffler function does not satisfy the composition property: \begin{equation}E_{\alpha}(x)E_{\alpha}(y)\neq E_{\alpha}(x+y),\end{equation} but it can be observed that the function \cite{16,17,18}: \begin{equation}E_{\alpha}( x^{\alpha})=\Sigma_{i=0}^{\infty}\frac{x^{\alpha i}}{\Gamma(\alpha i+1)}\end{equation}
does satisfy the composition property:
\begin{equation}E_{\alpha}( x^{\alpha})E_{\alpha}( y^{\alpha})=E_{\alpha}((x+y)^{\alpha}).\quad \alpha>0, \label{jumarie formulae}\end{equation}
Using of the above relation  (\ref{jumarie formulae}), so the equation (\ref{equadeC1}) becomes:
\begin{equation}D^{\alpha}_{t}C_{1}=2(1-\lambda)\gamma^{2} E_{\alpha,1}(\theta (2t)^{\alpha}).\label{equa2deC1}\end{equation}
Applying $J^{\alpha}$ on both sides of  equation (\ref{equa2deC1}), and using  integration of the Mittag-Leffler function relation \cite{1} (p. 25), we obtain:
\begin{equation}\frac{1}{\Gamma(\nu)}\int_{0}^{t}(t-\tau)^{\nu-1}E_{\alpha,\beta}(\eta t^{\alpha})t^{\beta-1}dt=t^{\beta+\nu-1}E_{\alpha,\beta+\nu}(\eta t^{\alpha}), \label{relation podolbny}\end{equation}
where $\alpha>0,\, \beta>0$ and $ \nu>0$, it leads  by taking $\alpha=\nu,\,\, \eta=\theta 2^{\alpha}$ and $ \beta=1$ to:
\begin{equation} J^{\alpha}E_{\alpha,1}(\theta (2t)^{\alpha})=t^{\alpha}E_{\alpha,\alpha+1}(\theta (2t)^{\alpha}).\label{relation podolbny2}\end{equation}According to  the following  relation, it yields: \begin{equation} J^{\alpha}D^{\alpha}_{t}C_{1}(t)=C_{1}(t)-C_{1}(0),\qquad 0<\alpha<1,\end{equation}we obtain
\begin{equation}C_{1}(t)=2(1-\lambda)\gamma^{2}t^{\alpha} E_{\alpha,\alpha+1}(\theta (2t)^{\alpha})+C_{1}(0).\label{equa3deC1}\end{equation}
We  assume $C_{1}(0)=0$. Hence, the obtained  solution of fractional ordinary differential system (\ref{system}) yields the following exact solution of  the nonlinear time fractional  modified Kuramoto-Sivashinsky  equation (\ref{mkseq}):
\begin{equation}u(t,x)=2(1-\lambda)\gamma^{2}t^{\alpha} E_{\alpha,\alpha+1}(\theta (2t)^{\alpha})+E_{\alpha,1}(\theta t^{\alpha})(\cos\gamma x+\sin\gamma x),\label{solutionmkseq}\end{equation}
where $\gamma=\sqrt{\frac{1-\lambda}{\lambda}}$  and $\theta=\gamma^{2}-\gamma^{4}$.\end{section}
\begin{section}{SOME PARTICULAR CASES}
In this section, we extract some particular cases, precisely exact solutions corresponding to $\lambda=\frac{1}{m},$ with $ m\in\mathbb{N^{\star}}-\{1\}$.\\

\noindent{}\textbf{Case 1\,\,$\lambda=\frac{1}{2}$}\\ This particular value of $\lambda$ leads to $\gamma=1$ and $\theta=0$. Consequently, an exact solution of the studied fractional equation (\ref{mkseq}) is given by:
\begin{equation}u(t,x)=\frac{t^{\alpha}}{\Gamma(\alpha+1)}+\cos x+\sin x,\,\, 0<\alpha\leq1.\label{solutparticliere1}\end{equation}

\noindent{}\textbf{Case 2\,\,$\lambda=\frac{1}{m},\, m>2$}\\

\noindent{}In this case we obtain that $\gamma=\sqrt{m-1}$ and $\theta=(m-1)(2-m).$ The constructed exact solution takes the form:
\begin{eqnarray}u_{\alpha,m}(t,x)&=&2^{1-\alpha}\frac{(m-1)}{m(2-m)}\left\{E_{\alpha,1}\left[(m-1)(2-m)2^{\alpha}t^{\alpha}\right]-1\right\}\nonumber\\
&+&E_{\alpha,1}\left[(m-1)(2-m)t^{\alpha}\right]\cos\sqrt{m-1}x\nonumber\\
&+&E_{\alpha,1}\left[(m-1)(2-m)t^{\alpha}\right]\sin\sqrt{m-1}x.
\label{solutparticliere2}\end{eqnarray}
\noindent{}Now we look for  solutions of nonlinear time fractional equation (\ref{mkseq}) corresponding to   $\alpha=1$ and $\alpha=\frac{1}{2}$.\\

\noindent{}\textbf{Subcase 2.1 $\lambda=\frac{1}{m},\, m>2, \alpha=1$.}\\

\noindent{}According to the relation \begin{equation}E_{1,1}(z)=E_{1}(z)=e^{z},\end{equation} the corresponding exact solution of equation (\ref{mkseq}) is obtained to be of the following form:
\begin{eqnarray}u_{1,m}(t,x)&=&\frac{m-1}{m(2-m)}\left\{e^{2(m-1)(2-m)t}-1\right\}+e^{(m-1)(2-m)t}\cos\sqrt{m-1}x\nonumber\\
&+&e^{(m-1)(2-m)t)}\sin\sqrt{m-1}x.
\label{solutparticliere21}\end{eqnarray}

\noindent{}\textbf{Subcase 2.2 $\lambda=\frac{1}{m},\, m>2, \,\alpha=\frac{1}{2}$.}\\

\noindent{}According to the relation:
\begin{equation}E_{\frac{1}{2},1}(z)=E_{\frac{1}{2}}(z)=e^{z^{2}}\left(1+erf(z)\frac{2}{\sqrt{\pi}}\int_{0}^{z}e^{-y^{2}}\,\,dy\right),\end{equation}
where \begin{equation}erf(z)=\frac{2}{\sqrt{\pi}}\int_{0}^{z}e^{-y^{2}}\,\,dy,\end{equation} the corresponding exact solution of equation (\ref{mkseq}) is obtained in this subcase  to be of the form:
\begin{eqnarray}u_{\frac{1}{2},m}(t,x)&=&\frac{\sqrt{2}(m-1)}{m(2-m)}\left\{E_{\frac{1}{2}}\left[\sqrt{2}(m-1)(2-m)\sqrt{t}\right]-1\right\}\nonumber\\
&+&e^{(m-1)(2-m)t}\left(\cos\sqrt{m-1}x+\sin\sqrt{m-1}x.\right).
\label{solutparticliere21}\end{eqnarray}
We  end this section by giving some corresponding graphs of some found particular solutions.
\begin{figure}[htbp]
  \begin{minipage}[b]{0.48\linewidth}
   \centering
   \includegraphics[width=6cm,height=5cm]{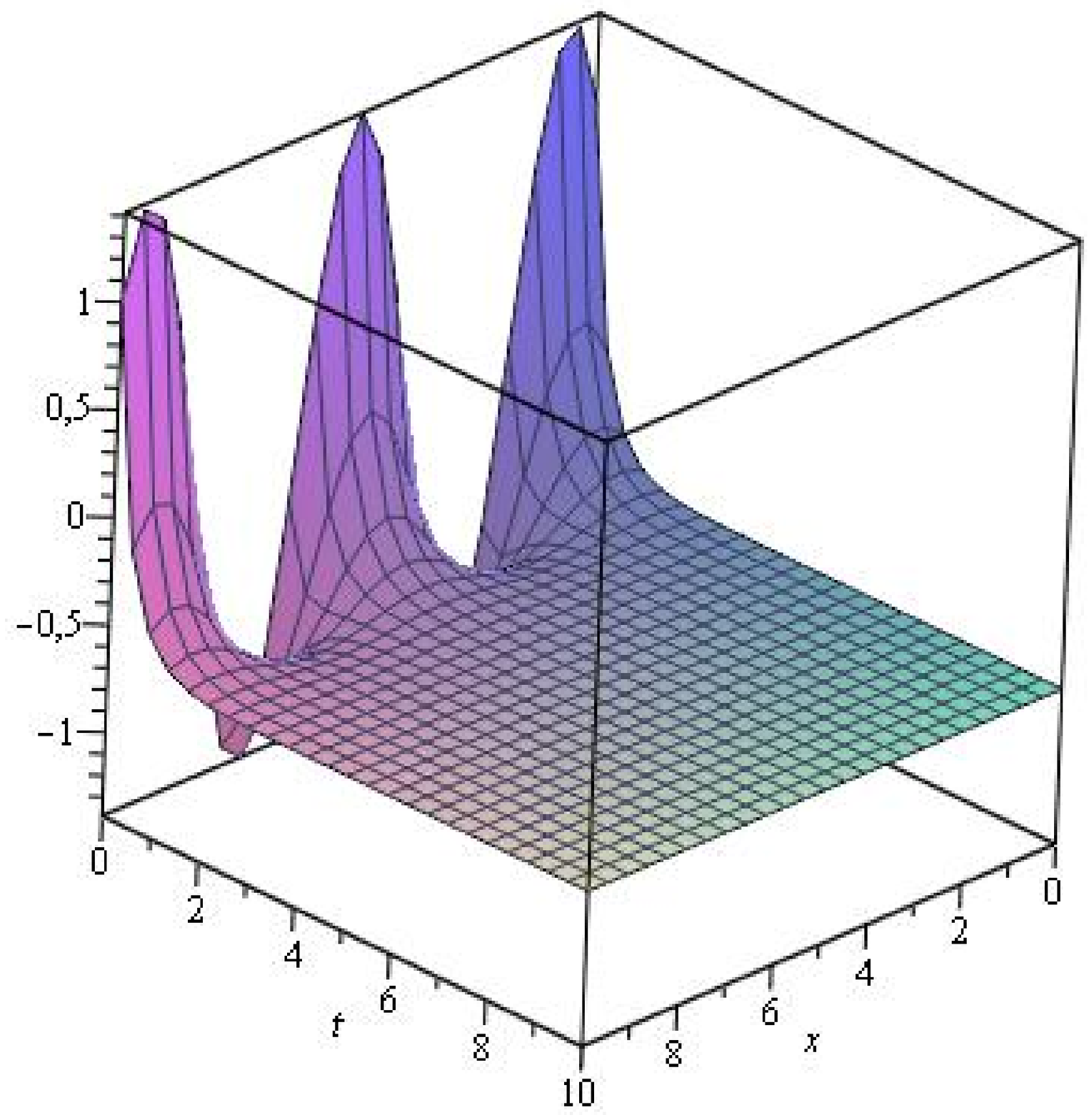}
   \caption[Fig.]{Solution with $\alpha=1$ and $m=3$}
  \end{minipage}
\hfill
\begin{minipage}[b]{0.48\linewidth}
   \centering
   \includegraphics[width=6cm,height=5cm]{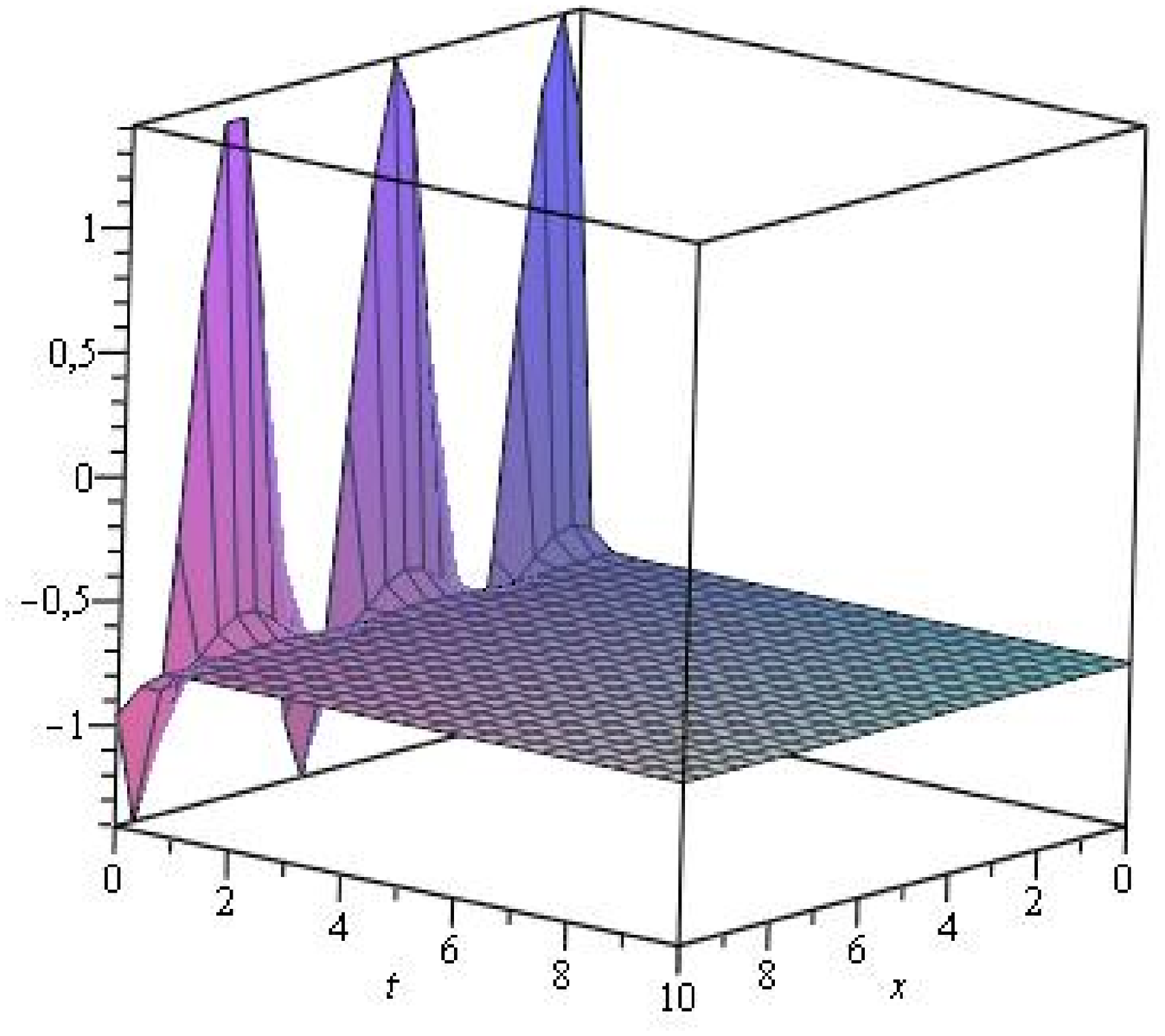}
   \caption{Solution with $\alpha=1$ and $m=4$}
  \end{minipage}
   \label{fig:ma_fig}
\end{figure}

\begin{figure}[htbp]
  \begin{minipage}[b]{0.48\linewidth}
   \centering
   \includegraphics[width=6cm,height=5cm]{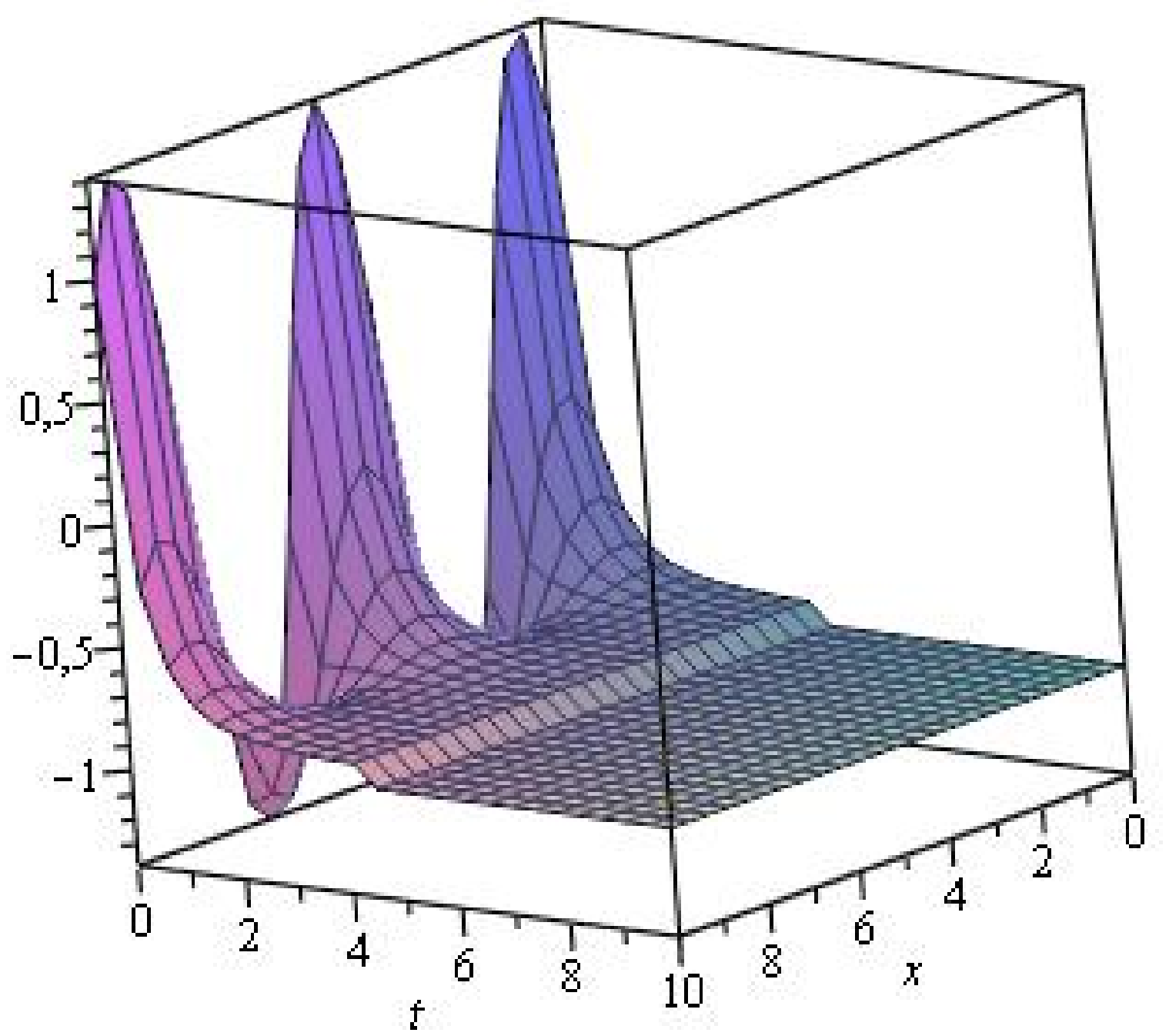}
   \caption[Fig.]{Solution with $\alpha=\frac{1}{2}$ and $m=3$}
  \end{minipage}
\hfill
\begin{minipage}[b]{0.48\linewidth}
   \centering
   \includegraphics[width=6cm,height=5cm]{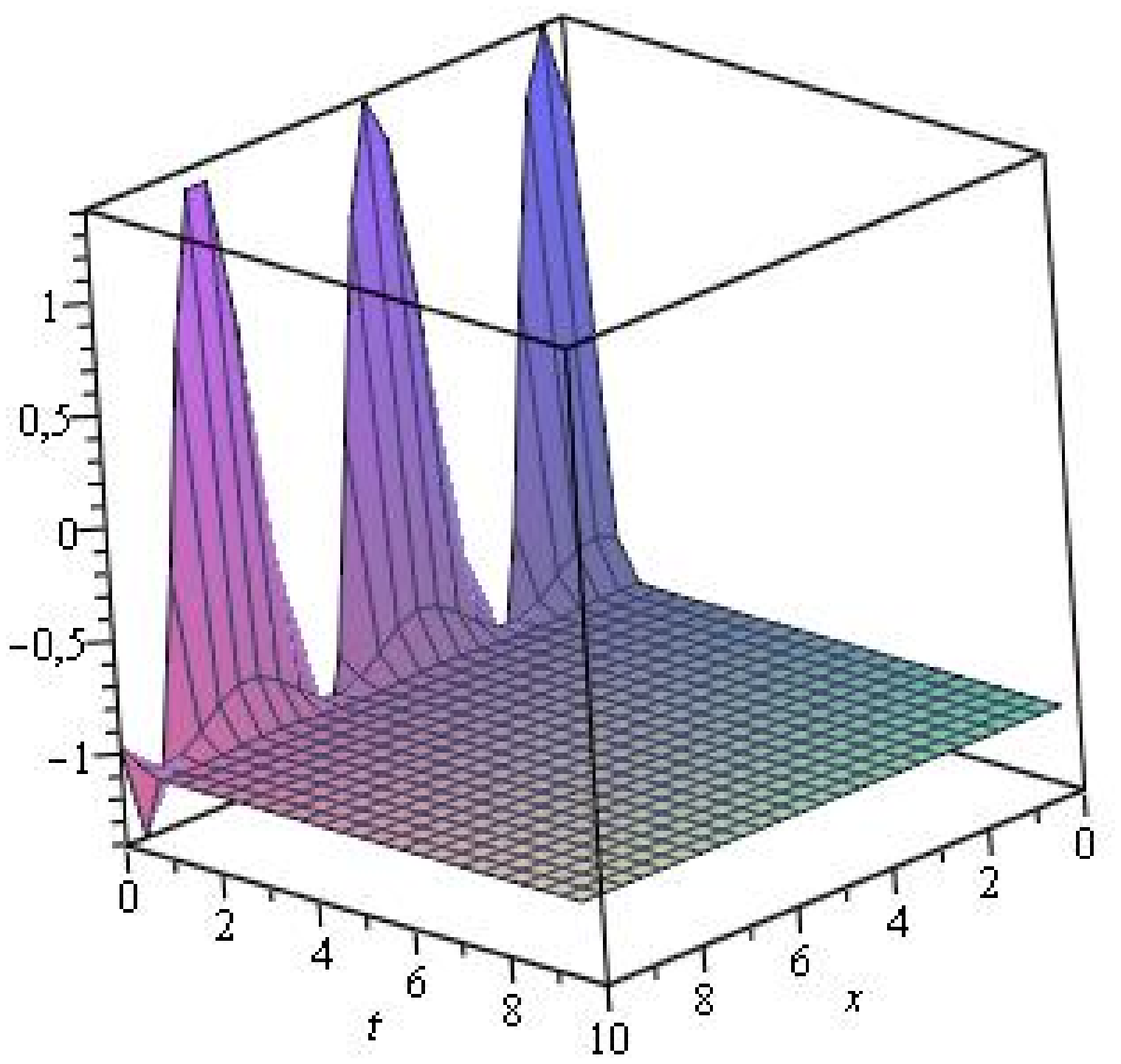}
   \caption{Solution with $\alpha=\frac{1}{2}$ and $m=4$}
  \end{minipage}
   \label{fig:ma_fig}
\end{figure}\vspace{-2cm}
\begin{figure}[htbp]
  \begin{minipage}[b]{0.48\linewidth}
   \centering
   \includegraphics[width=6cm,height=5cm]{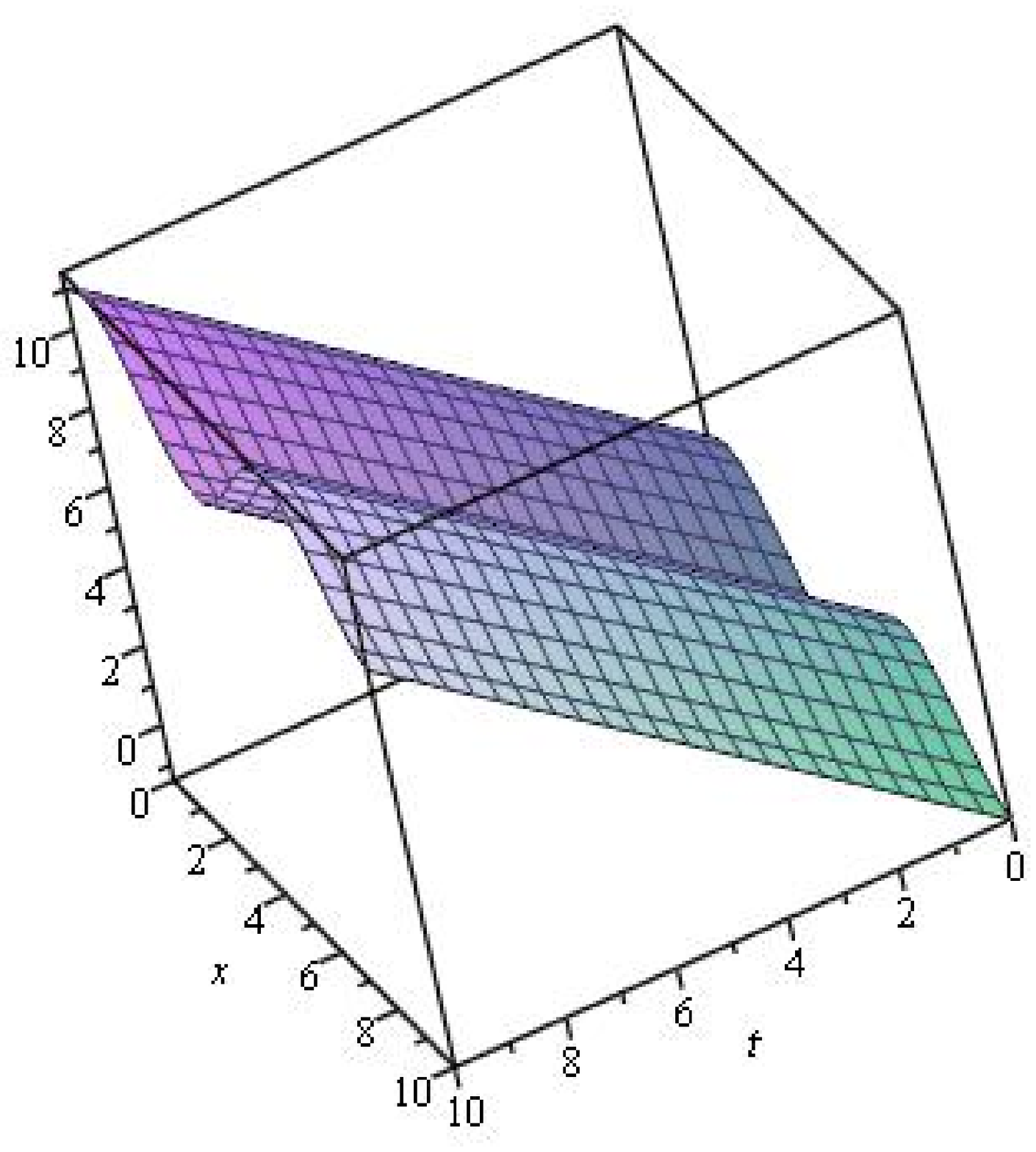}
   \caption[Fig.]{Solution with $\alpha=1$ and $\lambda=\frac{1}{2}$}
  \end{minipage}
\hfill
\begin{minipage}[b]{0.48\linewidth}
   \centering
   \includegraphics[width=6cm,height=5cm]{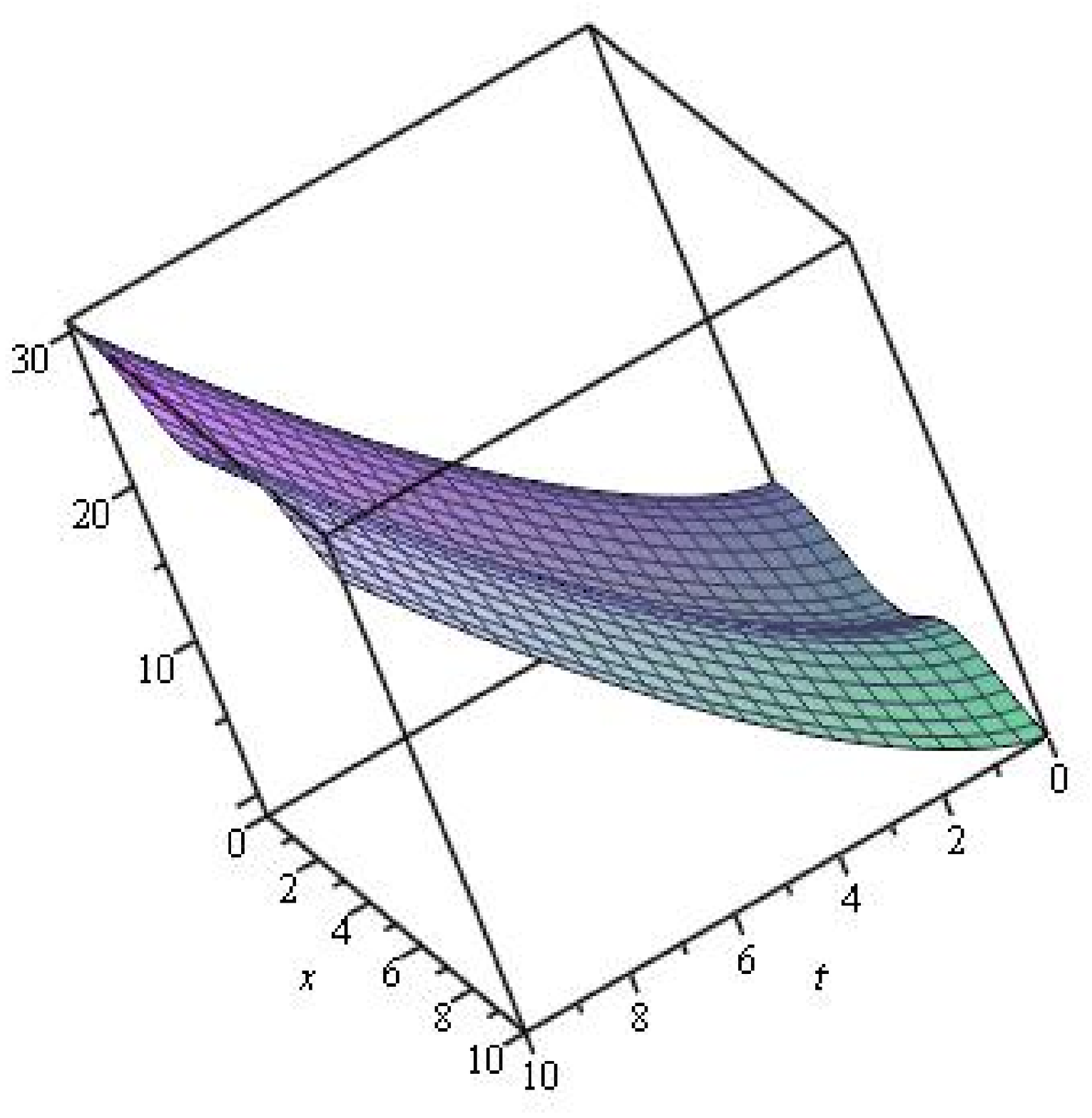}
   \caption{Solution with $\alpha=\frac{1}{2}$ and $\lambda=\frac{1}{2}$}
  \end{minipage}
   \label{fig:ma_fig}
\end{figure}
\vspace{-2cm}
\end{section}
\newpage
\begin{section}{Conclusion}
Here, by using the Laplace transform method and some basic properties of the Mittag-Leffler function, we succeed to solve the obtained reduced system of ordinary fractional equations. Consequently, the invariant subspace method was appropriate to construct some exact solutions of the time fractional nonlinear modified  Kuramoto-Sivashinsky equation. Finally, we note that the method used in this paper can be extended to obtain exact solutions of other nonlinear time fractional differential equations.
\end{section}

\end{document}